\documentclass[11pt]{article}
\usepackage[sectionbib,square, numbers, comma, sort&compress]{natbib} 
\usepackage{amsmath}
\usepackage{amssymb}
\usepackage{amsfonts}
\usepackage{bm}
\usepackage{amsthm}
\usepackage[T1]{fontenc}
\usepackage{mathrsfs}
\usepackage{pbox}
\usepackage{multirow}
\usepackage{caption}
\usepackage{comment}
\usepackage{xcolor}
\usepackage{paralist}
\usepackage{enumitem}
\usepackage{float}
\date{}
\allowdisplaybreaks
\begin{document}
\newcommand{\ot}{\otimes}
\newcommand{\m}{\mu}
\newcommand{\om}{\omega}
\newcommand{\ga}{\gamma}
\newcommand{\ze}{\zeta}
\newcommand{\N}{N_{\F_{2^{2t}}/\F_{2^2}}}
\newcommand{\Nn}{N_{\F_{2^{n}}/\F_{2^2}}}
\newcommand{\wt}{\operatorname{wt}}
\renewcommand{\l}{\lambda}
\renewcommand{\arraystretch}{1.3}
\newcommand{\F}{\mathbb{F}}
\newtheorem{fact}{Fact}

\title{On Some Permutation Binomials and Trinomials Over $\F_{2^n}$}
\author {{\bf Srimanta Bhattacharya}\\
Centre of Excellence in Cryptology,\\
Indian Statistical Institute,\\ Kolkata.\\
E-mail: mail.srimanta@gmail.com\\
\and
{\bf Sumanta Sarkar}\\
TCS Innovations Labs,\\Hyderabad.\\
E-mail: sumanta.sarkar@gmail.com\\
}
\maketitle
\newtheorem{theorem}{\bf Theorem}[section]
\newtheorem{proposition}[theorem]{\bf Proposition}
\newtheorem{lemma}[theorem]{\bf Lemma}
\newtheorem{corollary}[theorem]{\bf Corollary}
\newtheorem{definition}{\bf Definition}[section]
\newtheorem{observation}{Observation}
\theoremstyle{definition}
\newtheorem{example}{\bf Example}[section]
\theoremstyle{remark}
\newtheorem{remark}{Remark}[theorem]
\setlength{\abovedisplayskip}{0pt}
\setlength{\belowdisplayskip}{0pt}
\setlength{\abovedisplayshortskip}{0pt}
\setlength{\belowdisplayshortskip}{0pt}
\setlist[itemize]{noitemsep, topsep=0pt}
\setlist[description]{noitemsep, topsep=0pt}
\setlist[enumerate]{noitemsep, topsep=0pt}
\setlist{nolistsep}
\setlength{\bibsep}{0pt}
\setlength{\tabcolsep}{0.2em}
\begin{abstract}
In this work, we completely characterize 
\begin{inparaenum}[(\itshape i\upshape)]
\item permutation binomials of the form $x^{{{2^n -1}\over {2^t-1}}+1}+ ax \in \F_{2^n}[x], n = 2^st, a \in \F_{2^{2t}}^{*}$, and
\item permutation trinomials of the form $x^{2^s+1}+x^{2^{s-1}+1}+\alpha x \in \F_{2^t}[x]$, 
\end{inparaenum} where $s,t$ are positive integers.
The first result, which was our primary motivation, is a consequence of the second result. The second result may be of independent interest.
\end{abstract}
{\bf Keywords:}
Finite field, permutation binomial, permutation trinomial.
\section{Introduction}
\subsection{Motivation}
Let $\F_q$ be a finite field of $q$ elements.  A polynomial $f \in \F_q[x]$ is a {\em permutation polynomial} (PP) of $\F_q$ if the function $f: a \mapsto f(a), a \in \F_q$, is a permutation of $\F_q$. PPs represent purely combinatorial objects, namely permutations. This influences their algebraic properties making PPs theoretically appealing. Besides their theoretical importance, PPs have been considered in several practical contexts such as in cryptography, coding theory, combinatorial designs, etc. Motivated by their theoretical as well as practical significance, there seems to be a renewed interest in PPs in recent literature (see \cite{MP13,Hou15} for a recent account).\par
We term a class of PPs {\em characterized} if both necessary and sufficient conditions (to be PP) are available for polynomials belonging to that class. Characterization of PPs is an extremely important and challenging open question, and perhaps forms the crux of the research on PPs. Several classes of PPs have already been characterized; well-known characterized classes include linearized polynomials, Dickson polynomials, etc. However, despite considerable attention the problem remains unsolved in general.\par   
Among polynomials, monomials are easily characterized for their permutation properties; the monomial $x^d$ is a PP of $\F_q$ if and only if $\gcd(d,q-1) = 1$. However, the problem is already difficult for binomials of the form $x^d+ax$. In fact, precise characterization is not available even for binomials of the specific form $x^{{q-1 \over d}+1}+ax \in \F_q[x]$, where $d$ is a divisor of $q-1$. These binomials belong to the class of {\em cyclotomic mapping polynomials}\footnote{These are polynomials of the form $x^rf(x^{q-1\over d})$, and represent mappings of the factor group $\F_{q}^{*} / C_d$ to itself, where $C_d$ is the subgroup of $\F_{q}^{*}$ of index $d$ (see  \cite{Eva} for further details). We review relevant characterization results of this class of binomials in Section \ref{exres}.} and are relatively well studied with respect to their permutation properties starting with the work of Carlitz (\cite{Car62}), who first showed their existence for sufficiently large $q$ with respect to $d$. Also, these are very closely related to {\em complete mappings / orthomorphisms}. \footnote{A polynomial $f(x)\in \F_q[x]$ is called complete mapping if both $f(x)$ and $f(x)+x$ are PPs of $\F_q$, and orthomorphism if both $f(x)$ and $f(x)-x$ are PPs; for even characteristic both are same. Complete mappings / orthomorphisms are useful for construction of {\em mutually orthogonal latin squares} (see \cite{Eva,Zie13a}).}\par 
These considerations motivate us to address the problem of characterization of these {\em permutation binomials} (PBs) as the next non-trivial open case (after monomials). More specifically, in this work, we characterize PBs of the form $x^{{2^n-1\over 2^t-1}+1}+ax$ over $\F_{2^n}, n= 2^st, a \in \F_{2^{2t}}^{*}$, where $s,t$ are positive integers. As an additional result, we obtain complete characterization of permutation trinomials of the form $x^{2^s+1}+x^{2^{s-1}+1}+\alpha x \in \F_{2^t}[x]$; in fact, our first result on characterization of PBs follows from this result.\par 
  In the next subsection, we describe our contribution in more detail in the context of relevant existing results. Then we state the results that we require for our proofs. Finally, in Section \ref{proofs}, we give proofs of our results.\par

\subsection{Existing Results and Statement of Our Contribution}
\label{exres}
In this part, we discuss subclasses of PBs of the form $x^{{q-1 \over d}+1}+ax \in \F_q[x]$, $q$ even, that have already been characterized. Here, we point out that well-known Hermite-Dickson criteria (Theorem \ref{hdthm}) and Wan-Lidl criteria (Theorem \ref{permth1}), which we use as tools in our proofs, are themselves characterizations of larger classes of PPs, of which the PBs of the form $x^{{q-1 \over d}+1}+ax$ are a subclass. In fact, Hermite-Dickson criteria is a characterization for the class of all PPs. However, we are interested in characterizations which are more direct and explicit. Below we discuss cases where such explicit characterizations are available.\par
In \cite{Wan02}, the author characterized PBs of the form $x^r(x^{{q-1 \over d}s}+1) \in \F_{q}[x]$ for $d \in \{3,5\}$ and positive integers $r,s$. In \cite{AW05}, characterization for $d=7$ was obtained.\footnote{The case of $d=2$ was settled in \cite{NR82}. However, it is relevant for fields of odd characteristic.} \footnote{In \cite{Wan07, Wan10}, the author characterized these PBs for any $d$ in terms of Lucas sequences. However, as we have stated before, we are interested in more explicit characterization.}\par
In \cite{Zie13a}, the author characterized PBs of the form \begin{inparaenum}[(\itshape i \upshape)]\item $x^{{q^2-1 \over q-1}+1}+ax$ over $\F_{q^2}$, \item $x^{{q^3-1 \over q-1}+1}+ax$ over $\F_{q^3}$,\end{inparaenum}\footnote{PBs of the form $x^{2{q^2-1\over q-1}+1}+ax$ over $\F_{q^2}$ were also characterized in the same work.} for all characteristic. This characterization generalizes the results of \cite{GC15, TZH14, WLHZ14} pertaining to the PBs of this form. In \cite{Zie13a}, the author used a variant of Theorem \ref{permth1} to reduce the polynomials (in both the cases) into a polynomial of low degree (degree $\leq 5$) over $\F_q$; this reduction preserves permutation property. Then he used characterizations of such polynomials from Dickson's table.  The same approach was used in \cite{WLHZ14}. In \cite{BZ15}, the authors reproved the above results using different techniques; along with Wan-Lidl criteria (the version given in \cite{NR82}), they used criteria for solvability of bivariate equations over $\F_q$.\par
In another line of work, characterization of PBs of the form $x^r(x^{q^2-1\over q+1}+a)$ over $\F_{q^2}$ have been considered. In \cite{Zie13b}, the author characterized this class under the restriction that $a$ is a $(q+1)$-th root of unity. As for more specific classes, PBs of the form $x^{2(q-1)+1}+ax$ and $x^{3(q-1)+1}+ax$ over $\F_{q^2}$ have been characterized in \cite{Hou15b} and \cite{HL15} respectively. \par
\subsubsection*{Statement of our results}
In this context, our first result is explicit characterization of PBs of the form $x^{{{2^n -1}\over {2^t-1}}+1}+ ax$, $a \in \F_{2^{2t}}^{*},n = 2^st$. More precisely, our result is the following.
 \begin{theorem}
\label{binneg}
Let $s, t$ be positive integers and $n = 2^st$. 
Then the polynomial $x^{{{2^n -1}\over {2^t-1}}+1}+ ax$, $a \in \F_{2^{2t}}^{*}$, is a PP of $\F_{2^n}$ if and only if 
\begin{inparaenum}[{\upshape (}i{\upshape )}]\item $t$ is odd,\item $s\in\{1, 2\}$, and \item $a \in \om\F_{2^t}^{*} \cup \om^2 \F_{2^t}^*$, where $\om \in \F_{2^2}$ is a root of the equation $\om^2+\om+1=0$.
\end{inparaenum}
\end{theorem}
Theorem \ref{binneg} immediately leads to the following corollary.
\begin{corollary}
\label{chkcor}
Let $n= 2^st, s \in \{1, 2\}$, and $t$ be odd, then the number of $a \in \F_{2^{2t}}$ such that $x^{{{2^n -1}\over {2^t-1}}+1}+ ax \in \F_{2^{2t}}[x]$ is a PB is $2(2^t-1)$. 
\end{corollary}
Here, we highlight the fact that while we restrict $a$ to the subfield $\F_{2^{2t}}$ of $\F_{2^n}$, the setting for $n (= 2^st)$ is considerably general than the previously discussed cases of \cite{Zie13a, WLHZ14, TZH14, BZ15}; in all these cases $n$ is of the form $2^{st}$ for specific values of $s$.\par
The case $s=1$ of Theorem \ref{binneg} was proven earlier in \cite{CK08,SBC12}. Here, we mention that the main thrust of Theorem \ref{binneg} came from the recent work \cite{WLHZ14}, where PBs of the form $x^{{2^{st}-1 \over 2^t-1}+1}+ax \in\F_{2^{st}[x]}$ were investigated for $s \in \{3,4,6,10\}$ under certain restrictions on $t$ (see also \cite{TZH14}); there, the authors obtained sufficient conditions for the above classes of PBs. In particular, from their results, we observed that while there are $a \in \F_{2^{2t}}^{*}$ such that $x^{{2^{4t}-1 \over 2^t-1}+1}+ax$ is a PB of $\F_{2^{4t}}$, there are no such $a \in \F_{2^{2t}}^{*}$ such that $x^{{2^{8t}-1 \over 2^t-1}+1}+ax$ is a PB of $\F_{2^{8t}}$. Theorem \ref{binneg} is a generaliztion of this observation.  \par
In a simultaneous and independent work (\cite{BZ15a}), the authors have characterized PBs of the form $x^{{2^{4t}-1 \over 2^t-1}+1}+ax$ over $\F_{2^{4t}}$. They have shown that for $t(\geq 4)$ even, there does not exist any PB of this form; for $t(\geq 3)$ odd, they have characterized all PBs of this form. So, there is overlap of this result with Theorem \ref{binneg} for the case $s=2, a\in \F_{2^{2t}}$. However, their approach in this case is similar to \cite{BZ15}, and is different from ours.\par
 Our approach (discussed in more detail in Section \ref{proofs}) in the proof of Theorem \ref{binneg} is different from those of \cite{Zie13a,WLHZ14, BZ15, BZ15a} to some extent. In our case, we use Wan-Lidl criteria to reduce (preserving premutation property) the polynomial $x^{{{2^{n} -1}\over {2^t-1}}+1}+ ax\in \F_{n}[x]$, $n= 2^st, a \in \F_{2^{2t}}^{*}$, to a trinomial of the form $x^{2^s+1}+x^{2^{s-1}+1}+\alpha x \in \F_{2^t}$. To characterize permutation trinomials of the form $x^{2^s+1}+x^{2^{s-1}+1}+\alpha x$  we use Hermite-Dickson criteria (Theorem \ref{hdthm}) along with known characterization of low degree PPs from Dickson's table. More formally, we have the following result for these permutation trinomials. \par  
\begin{theorem}
\label{trith}
Let $s, t$ be positive integers. Then the polynomial $x^{2^s+1}+x^{2^{s-1}+1}+\alpha x \in \F_{2^t}[x]$ is a PP of $\F_{2^t}$ if and only if \begin{inparaenum}[{\upshape (}i{\upshape )}]\item $t$ is odd,\item $\alpha = 1$, and \item $s \in \{1, 2\}$.
\end{inparaenum} 
\end{theorem}
Our main motivation was to characterize PBs of the form $x^{{2^n-1 \over 2^t-1}+1}+ax \in \F_{2^n}[x]$ which led to Theorem \ref{binneg} and we use Theorem \ref{trith} for this purpose. However, Theorem \ref{trith} may be of independent interest, especially since permutation properties of trinomials are much less known (see \cite{DQWYY15, Hou15a}), and here we only hope that this result, due to the intersting form of the trinomial, proves to be a valuable addition to the set of few existing results on permutation trinomials.\par  
\subsection{Useful Results}
Hermite-Dickson criteria and Lucas' theorem (Theorem \ref{lucas}) are the main tools in our proof of Theorem \ref{trith}; in fact, we use a corollary (Corollary \ref{lucascor}) of Lucas' theorem. Finally, we derive Theorem \ref{binneg} from Theorem \ref{trith} using Wan-Lidl criteria.\par
\subsubsection*{Hermite-Dickson criteria:}
The first non-trivial characterization for the class of all PPs is given by the following result, commonly known as Hermite-Dickson criteria.
\begin{theorem}[see \cite{LN97}]
\label{hdthm}
A polynomial $f \in \F_q[x]$ is a PP if and only if \linebreak
\begin{inparaenum}[\upshape{(}\itshape i \upshape)]
\item $f$ has exactly one root in $\F_q$,
\item $f^t \mod (x^q-x)$ has degree less than $q-1$, for $1 \leq t \leq q-2, p\nmid t$, where $p$ is the characteristic of $\F_q$.
\end{inparaenum}
\end{theorem}
\begin{remark}
In the above theorem, one can remove the condition $p \nmid t$; we will do so in the proof of Theorem \ref{trith}.
\end{remark}
\subsubsection*{Lucas' theorem for multinomials:}
\begin{theorem}[Lucas (see \cite{LN97})]
\label{lucas}
Let $p$ be a prime, and $n, r_1, r_2, \ldots, r_t$ be nonnegative integers such that
\begin{align*}
& n = d_0+d_1p+d_2p^2+\ldots+d_sp^s ~ (0\leq d_i\leq p-1, \forall~ 0 \leq i\leq s),\\
&r_j = d_{j0}+d_{j1}p+d_{j2}p^2+\ldots+d_{js}p^s ~ (0\leq d_{ji} \leq p-1, \forall~ 1 \leq j \leq t, \forall~ 0\leq i \leq s).
\end{align*}
Then
\vspace{5pt}
$$
\binom{n}{r_1, r_2, \ldots, r_t} = \binom{d_0}{d_{10}, d_{20}, \ldots, d_{t0}} \cdots \binom{d_s}{d_{1s}, d_{2s}, \ldots, d_{ts}} \mod (p).
$$
\end{theorem}
\vspace{1pt}
\begin{corollary}
\label{lucascor} 
With the notation of Theorem \ref{lucas}, it follows that\linebreak $\binom{n}{r_1, r_2, \ldots, r_t} \not = 0 \mod (p)$ iff $\sum_{i=1}^{t}d_{ij} = d_j,  \forall~ 0 \leq j \leq s$.
\end{corollary}
\subsubsection*{Wan-Lidl criteria:}
The following theorem from \cite{LW91} is a fundamental tool for analyzing permutation properties of polynomials of the form $x^rf(x^{q-1\over d})$. It is quite useful in that it reduces permutation property of the polynomial $x^rf(x^{q-1\over d}) \in \F_q[x]$ into permutation property of a related polynomial over a smaller subset, the set of $d$-th roots of unity, of $\F_q$. In some cases, choice of the smaller subset to be the multiplicative group of a subfield reduces the original polynomial ($x^rf(x^{q-1 \over d})$) to a polynomial whose permutation properties are known.\footnote{This approach was taken in \cite{Zie13a, WLHZ14} (see also \cite{AGW11} and references therein) .} The theorem was reproven at various other places. Here, it should be noted that the case $r=1$ was proven earlier in \cite{NR82}; in fact, for our purpose the case $r=1$ is sufficient.
\begin{theorem}[\cite{LW91}]
\label{permth1}
Let $d$ and $r$ be two positive integers and $q$ be a prime power such that $d$ divides $q-1$. Let
$\ga$ be a primitive element in $\mathbb{F}_q$ and assume $f$
is a polynomial in $\mathbb{F}_q[x]$. Then
$g(x) = x^rf(x^{\frac{q-1}{d}})$
is a PP of $\mathbb{F}_q$ if and only if the
following conditions are satisfied.
\begin{inparaenum}[(a)]
\item $\gcd(r, \frac{q-1}{d}) = 1$,
\item \mbox { for all } $i$,
$0 \leq i < d$, $f(\ga^{i\frac{q-1}{d}}) \neq 0$,
\item \mbox { for all } $j$, $0 \leq i < j < d$,
$g(\ga^i)^{\frac{q-1}{d}} \neq g(\ga^j)^{\frac{q-1}{d}}$.
\end{inparaenum}
\footnote{Conditions (b) and (c) can be written together (see \cite{Zie13b}) as the condition: $x^rf(x)^{q-1\over d}$ permutes the set $\{a \in \F_q: a^d=1\}.$}
\end{theorem}
\par
\section{Proofs}
\label{proofs}
Let $a = \sum_{i=0}^{\ell}a_i2^i$ be the $2$-adic representation of $a$, then we denote by $a_i$ the $i$-th binary digit (bit) of $a$. Also, let $\wt(a) = \lvert \{i \vert a_i \neq 0\}\rvert$.\par
\subsection{Proof of Theorem \ref{trith}}
{\em Proof.} First, we note that for positive integers $s > t$, the polynomial \linebreak $f(x)=x^{2^s+1}+x^{2^{s-1}+1}+\alpha x \in \F_{2^t}[x]$ can be reduced modulo $x^{2^t}+x$ to a polynomial $x^{2^{s'}+1}+x^{2^{s'-1}+1}+\alpha x$, with $s' < t$, which induces identical mapping on $\F_{2^t}$. Hence, it is sufficient to consider the cases with $s \leq t$.\par
Next, we consider the cases corresponding to $s\in \{1,2, t\}$, and for these cases we directly refer to the work of Dickson \cite{Dic96} (see also \cite{LN97}), where all PPs of degree $\leq 5$ for all characteristics were characterized.\par 
For $s=t$, we have $f(x) = x^2+x^{2^{t-1}+1}+\alpha x \bmod (x^{2^t}+x)$. Note that $f(x)$ is a PP of $\F_{2^t}$ if and only if $g(x)= f(x)^2 \bmod (x^{2^t}+x)$ is a PP of $\F_{2^t}$. Now, $g(x)=f(x)^2 \mod (x^{2^t}+x) = x^4+x^3+\alpha^2x^2$. Characterization of \cite{Dic96}, along with direct calculations for the cases $t \in \{1,2\}$, implies that $g(x)$ is a PP of $\F_{2^t}$ iff $t=1$ and $\alpha =1$. \par
Similarly, it also follows from \cite{Dic96} that for the cases $s = 1$ and $2$, $f(x) = x^{2^s+1}+x^{2^{s-1}+1}+\alpha x \in \F_{2^t}[x]$ is a PP of $\F_{2^t}$ if and only if $t$ is odd and $\alpha = 1$. In fact, for $s=2, \alpha=1$, $f(x) = x^5+x^3+\alpha x$ is the Dickson polynomial $D_5(x,1)$ (see \cite{LN97}), which is a PP of $\F_{2^t}$ if and only if $\gcd(2^{2t}-1, 5) = 1$, which is true if and only if $t$ is odd.\par 
Now, we show that $f$ is not a PP for $t>s \geq 3$ by applying the Hermite-Dickson criteria (Theorem \ref{hdthm}). For this, first we raise $f$ to $2^t-3$ modulo $x^{2^t}+x$, and then to $2^t-4$  modulo $x^{2^t}+x$; we show that the degree of the resulting polynomial is $2^t-1$ in at least one of these two cases. Here, it is important to note that for any polynomial $g \in \F_{2^t}[x]$, exactly those terms whose exponents are multiples of $2^t-1$ reduce to the term with exponent $2^t-1$ when $g$ is reduced modulo $x^{2^t}+x$. More precisely, and specifically for our case, we note the following fact which will be used later.
\begin{fact}
Let $g = \sum_i a_ix^i \in \F_{2^t}[x]$, and let $g \bmod (x^{2^t}+x) = b_{2^t-1} x^{2^t-1}+\sum_{i=0}^{2^t-2}b_i x^i$. Then $b_{2^t-1} = \sum_{j = 0 \mod (2^t-1)}a_j$.
\end{fact}
Hence, we will be done if we can show that sum of the coefficients of the terms whose exponents are multiples of $2^t-1$ in the exapnsion of $f^{2^t-3}$ or $f^{2^t-4}$ is non-zero. For this, we first consider the expansion of $f^{2^t-3}$ and then of $f^{2^t-4}$; we show that if in the first case the sum is zero then it is non-zero in the second case. Though the approaches are similar in these two cases, they are not exactly same.\\
{\bfseries \textsf{Case 1. ${f^{2^t-3}}:$}} First, note that coefficient of a term whose exponent is $\ell(2^t-1), \ell\geq 1$, in the expansion of $f^{2^t-3}$  is $(\binom{2^t-3}{u,v,w} \mod(2)) \alpha^w$, where $0 \leq u, v, w\leq 2^t-3$ are such that the following conditions hold:
\begin{align}
u+v+w &= 2^t-3 \label{eq1},\\[3pt]
(2^s+1)u+(2^{s-1}+1)v+w &= \ell(2^t-1).\label{eq2}
\end{align}
Let $\mathcal{S} = \{(u,v,w,\ell) \vert u,v,w,\ell~\mbox{non-negative, and satisfies (\ref{eq1}) and (\ref{eq2})}\}$. Our goal is to find expression of the sum $\sum_{(u,v,w,\ell) \in \mathcal{S}} (\binom{2^t-3}{u,v,w}\mod (2)) \alpha^w$. For this, we split the sum into parts according to the value of $\ell$, and investigate contribution from  each part.\par 
Henceforth, for this case, whenever we write $\binom{2^t-3}{u,v,w}$, we implicitly assume values of $u,v,w$ satisfying, possibly along with some other constraints, (\ref{eq1}) and (\ref{eq2}) for some $\ell$ whose value will be clear from the context. Also, we have the following observation.
\begin{observation}
$1$st bit of $2^t-3$ is zero. Hence, if any of $u, v, w \in \{2, 3\} \mod (4)$ then following Corollary \ref{lucascor} $\binom{2^t-3}{u, v, w} = 0 \mod (2)$.
\label{obs3}
\end{observation}
Now, (\ref{eq2})$-$(\ref{eq1}) yields 
\begin{eqnarray}
2^su+2^{s-1}v = 2^t(\ell-1)-(\ell-3). \label{eq3}
\end{eqnarray}
Clearly, both $u$ and $v$ can not be zero at the same time, and since $t > s$, we have from (\ref{eq3}), $\ell = 3 \mod (2^{s-1})$. Also, from (\ref{eq1}) and (\ref{eq2}), $\ell \leq 2^s+1$. So, possible values of $\ell$ are $3$ and $2^{s-1}+3$. We consider the following two subcases based on these two values of $\ell$.
\begin{description}
\item[Subcase 1.1.] ${\ell = 3:}$ In this case, (\ref{eq3}) yields $v = 2^{t-s+2}-2u$. Depending on $\wt(u)$ we consider the following  subsubcases.
\begin{description}
\item[Subsubcase 1.1.1.] {${\wt(u) > 1}$:} Let $u= \sum_{j=i_1-k+1}^{i_1} 2^j+ \sum_{j=0}^{i_2}u_j2^j$, where $k \geq 1, i_1 \leq t-s, i_2 \leq i_1-k-1, u_j \in \{0, 1\} ~\text{for}~0 \leq j \leq i_2$, and if $k=1$ then at least one $u_j$ is non-zero (since $\wt(u) > 1$).\\ 
 So, $v = 2^{t-s+2} - \sum_{j=i_1-k+2}^{i_1+1} 2^j- \sum_{j=1}^{i_2+1}u_{j-1}2^j$. Hence, $v = 2^{i_1+2} - \sum_{j=i_1-k+2}^{i_1+1} 2^j- \sum_{j=1}^{i_2+1}u_{j-1}2^j \mod (2^{i_1+2})$, since $i_1+2 \leq t-s+2$. Now, $2^{i_1+2} - \sum_{j=i_1-k+2}^{i_1+1} 2^j = 2^{i_1-k+2}$, and $\sum_{j=1}^{i_2+1}u_{j-1}2^j < 2^{i_2+2} \leq 2^{i_1-k+1}$. So, $v \mod (2^{i_1+2}) \leq 2^{i_1-k+2} $, and $v \mod (2^{i_1+2})\\ > 2^{i_1-k+2} - 2^{i_1-k+1} = 2^{i_1-k+1}$. Again, we have the following two possibilities.
\begin{enumerate}[label=(\textbf{\itshape\roman*})]
\item{$v \mod (2^{i_1+2}) < 2^{i_1-k+2}$:} In this case, $(i_1-k+1)$-th bit of $v$ is $1$, since $v \mod (2^{i_1+2}) > 2^{i_1-k+1} $. So, $(i_1-k+1)$-th bits of both $u$ and $v$ are $1$. Hence, following Corollary \ref{lucascor}, $\binom{2^t-3}{u,v,w}=0 \mod(2)$.
\item {$v = 2^{i_1-k+2} \mod(2^{i_1+2})$:} For this case, we observe that $k > 1$. Since otherwise, at least one $u_j$ in the sum $\sum_{j=0}^{i_2}u_j 2^j$, appearing in the binary representation of $u$, is non-zero, which implies $v \neq 2^{i_1-k+2} \mod(2^{i_1+2})$, a contradiction. Now, for $k> 1$, $(i_1-k+2)$-th bit of both $u$ and $v$ are $1$. So, again $\binom{2^t-3}{u,v,w}=0 \mod (2)$. 
\end{enumerate}
\item[Subsubcase 1.1.2.] ${\wt(u) \leq 1}$: For $\ell=3$, (\ref{eq3}) implies $u \leq 2^{t-s+1}$. Also, if $u = 1$ then $v = 2^{t-s+2}-2u$, i.e., $v = 2 \mod (4)$. So, by Observation \ref{obs3}, $\binom{2^t-3}{u, v, w} = 0 \mod (2)$. For the remaining possible values of $u$, i.e., for $u = 0$ or $2^i$, with $2 \leq i \leq t-s+1$, we examine the bit patterns of $u, v, w$ in Table \ref{tab1}. For better understanding, we illustrate the case $t = 9, s = 3, i =4$ in Table \ref{tab2}.\par
\begin{table}[H]
\caption{Bit patterns for Subsubcase 1.1.2}
\label{tab1}
\hspace{2cm}
\begin{tabular}{|c|c|c|}
\hline
Values & \multicolumn{2}{c|}{Bit positions with $1$}\\
\hline
\multirow{3}{*}{\pbox{20 cm}{$u =0$,\\ $v = 2^{t-s+2}$,\\ $w =2^t - 2^{t-s+2} -3 $}} & $u$ & $\emptyset$\\ \cline{2-3} & $v$ & $\{t-s+2\}$\\ \cline{2-3} & $w$ & $\begin{aligned}\{r \vert &r = 0,\\
				      			      &2 \leq r \leq t-s+1,\\
							      &t-s+3 \leq r \leq t-1 \}\end{aligned}$ \\
\hline
\multirow{3}{*}{
\pbox{20 cm}{
$u =2^i$,\\ $v = 2^{t-s+2} - 2^{i+1}$,\\ $w = 2^t - 2^{t-s+2}+2^i -3$\\
  $(2 \leq i \leq t-s+1)$ }} & $u$ & $\{i\}$\\ 
\cline{2-3}
& $v$ & $\{r \vert i+1 \leq r \leq t-s+1\}$\\
\cline{2-3}
& $w$ & $\begin{aligned}\{r \vert &r=0, \\
				&2 \leq r \leq i-1,\\
				&t-s+2 \leq r \leq t-1\}\end{aligned}$\\
\hline
\end{tabular}
\end{table}
{
\setlength{\tabcolsep}{0.4em}
\begin{table}[htb]
\caption{Example bit patterns for Subsubcase 1.1.2}
\label{tab2}
\hspace{2 cm}
\begin{tabular}{|c|cccccccccc|}
\hline
Values & \multicolumn{10}{c|}{Bit representation}\\
\hline
$2^t-3 = 2^9-3$&$1$&$1$&$1$&$1$&$1$&$1$&$1$&$1$&$0$&$1$\\ 
$u = 2^4$&$0$&$0$&$0$&$0$&$0$&$1$&$0$&$0$&$0$&$0$\\ 
$v = 2^8-2^5$&$0$&$0$&$1$&$1$&$1$&$0$&$0$&$0$&$0$&$0$\\ 
$w = 2^9-2^8+2^4-3$&$0$&$1$&$0$&$0$&$0$&$0$&$1$&$1$&$0$&$1$\\ 
\hline
\end{tabular}
\end{table}
}
\par
From Table \ref{tab1}, it can be observed that for these $t-s+1$ values of $u$, none of $u, v, w$ has $1$ in their $1$st bit position, each of $u, v, w$ has $0$ in their $t$-th bit position, and exactly one among $u, v, w$ has $1$ in each of the remaining $t-1$ bit positions. Hence, it follows from Corollary \ref{lucascor} that for each of these $t-s+1$ values of $u$, $\binom{2^t-3}{u, v, w} = 1 \mod (2)$.  
\end{description} 
So, coefficient of the term with exponent $3(2^t-1)$ in the expansion of $f^{2^t-3}$ is
$\alpha^{2^t-2^{t-s+2}-3}(1+\sum\limits_{i=2}^{t-s+1}\alpha^{2^i})$.
\item[Subcase 1.2.] ${\ell = 2^{s-1}+3:}$ 
 For $\ell = 2^{s-1}+3$, (\ref{eq3}) yields $v = 2^{t-s+2}(2^{s-2}+1)-2u-1$. Hence, $v = 3 \mod (4)$ when $u = 0 \mod(2)$ (since $s \geq 3$ and $t > s$). Therefore, by Observation \ref{obs3}, $\binom{2^t-3}{u, v, w} = 0 \mod (2)$ when $u = 0 \mod(2)$. Next, from (\ref{eq1}) and (\ref{eq2}), we get $w = u - 2^{t-s+2} -2$. Hence, for $u = 1 \mod (4), w = 3 \mod (4)$, which, by Observation \ref{obs3}, implies $\binom{2^t-3}{u, v,w} = 0 \mod (2)$. Again using Observation \ref{obs3}, $\binom{2^t-3}{u, v,w} = 0 \mod (2)$ for $u = 3 \mod(4)$.
\end{description}
So, coefficient of the term with exponent $2^t-1$ in the expansion of $f^{2^t-3} \mod x^{2^t}+x$ is $\alpha^{2^t-2^{t-s+2}-3}(1+\sum\limits_{i=2}^{t-s+1}\alpha^{2^i})$. Hence, if $1+\sum\limits_{i=2}^{t-s+1}\alpha^{2^i} \neq 0$ then $x^{2^s+1}+x^{2^{s-1}+1}+x, 3 \leq s < t,$ is not a PP of $\F_{2^t}$. Otherwise, i.e., if
\begin{eqnarray}
\label{2k3condn}
\sum\limits_{i=2}^{t-s+1}\alpha^{2^i} = 1,
\end{eqnarray}
we consider the next case.\par
{\bfseries \textsf{Case 2. ${f^{2^t-4}}:$}} Similar to equations (\ref{eq1}), (\ref{eq2}), and (\ref{eq3}) from the previous case, we get from the expansion of $f^{2^t-4}$ the following set of equations:
\begin{align}
\setlength{\jot}{50pt}
u+v+w &= 2^t-4, \label{eq4}\\[5pt]
(2^s+1)u+(2^{s-1}+1)v+w &= \ell(2^t-1), \label{eq5}\\[5pt]
2^su+2^{s-1}v &= 2^t(\ell-1)-(\ell-4). \label{eq6}
\end{align}
As in the previous case, when we write $\binom{2^t-4}{u,v,w}$, we mean values of $u, v,w$ that satisfy (\ref{eq4}), (\ref{eq5}) (and thereby (\ref{eq6})) for some $\ell$ which is clear from the context. Similar to Observation \ref{obs3}, we have the following observation in this case.
\begin{observation}
$0$-th bit and $1$st bit of $2^t-4$ are zero. So, if any of $u, v, w \in \{1,2, 3\} \mod (4)$ then following Corollary \ref{lucascor} $\binom{2^t-4}{u,v,w} = 0 \mod (2)$.
\label{obs4}
\end{observation}
Next, following similar considerations as in Case 1, from (\ref{eq6}), we get $\ell \in \{4, 2^{s-1}+4\}$. Now, for $\ell = 2^{s-1}+4$, $v = 2^{t-s+1}(2^{s-1}+3)-2u-1$. Since, $t > s$, we have $v \in \{1, 3\} \mod(4)$. This implies, by Observation \ref{obs4}, $\binom{2^t-4}{u, v, w} = 0 \mod (2)$. So, we are left with the $\ell = 4$ case. Now, we consider its following subcases.
\begin{description}
\item[Subcase 2.1.] ${\wt(u)\leq 1:}$ In this case, $u = 0$, or $u = 2^i, 0 \leq i \leq t-s+1$ (upper bound on $i$ follows from (\ref{eq6})). Now, if $i \in\{0, 1\}$ then following Observation \ref{obs4}, $\binom{2^t-4}{u, v, w} = 0 \mod (2)$. Also, for $i = t-s+1$, $v = 2^{t-s+1}$. So, both $u$ and $v$ have $1$ in their $(t-s+1)$-th bit position, which again implies $\binom{2^t-4}{u, v, w} = 0 \mod (2)$ for $i = t-s+1$. For $u=0$ and for the remaining values of $i$, i.e., for $2 \leq i \leq t-s$, we show the bit patterns of $u, v,w$ in Table \ref{tab3}. For better understanding, we illustrate the case for $t = 11, s = 4, i =5$ in Table \ref{tab4}.\par\pagebreak
\begin{flushleft}
\begin{table}[h]
\caption{Bit patterns for Subcase 2.1}
\label{tab3}
\hspace{1.5cm}
\begin{tabular}{|c|c|c|}
\hline
Values & \multicolumn{2}{c|}{Bit positions with $1$}\\
\hline
\multirow{3}{*}{\pbox{20 cm}{$u =0$,\\ $v = 2^{t-s+2}+ 2^{t-s+1}$,\\ $w =2^t - 2^{t-s+2}-2^{t-s+1} -4 $}} & $u$ & $\emptyset$\\ \cline{2-3} & $v$ & $\{t-s+2, t-s+1\}$\\ \cline{2-3} & $w$ & $\begin{aligned}\big\{r \vert & 2 \leq r \leq t-s,\\
&t-s+3 \leq r \leq t-1 \big\} \end{aligned}$ \\
\hline
\multirow{3}{*}{
\pbox{20 cm}{
$u =2^i$,\\ $v = 2^{t-s+2}+2^{t-s+1} - 2^{i+1}$,\\ $w = 2^t - 2^{t-s+2}- 2^{t-s+1}+2^i -4$\\
  $(2 \leq i \leq t-s)$ }} & $u$ & $\{i\}$\\ 
\cline{2-3}
& $v$ & $\begin{aligned}\big\{r \vert & r = t-s+2, \\
& i+1 \leq r \leq t-s\big \}\end{aligned}$\\
\cline{2-3}
& $w$ & $\begin{aligned}\big\{r \vert & 2 \leq r\leq i-1,\\
&  r = t-s+1,\\
& t-s+3 \leq r \leq t-1\big\} \end{aligned}$\\
\hline
\end{tabular}
\end{table}
\end{flushleft}
\vspace{-.4cm}
{
\setlength{\tabcolsep}{0.4em}
\begin{table}[h]
\captionsetup{justification=raggedright}
\caption{Example bit patterns for Subcase 2.1}
\label{tab4}
\hspace{1.5 cm}
\begin{tabular}{|c|cccccccccccc|}
\hline
Values & \multicolumn{12}{c|}{Bit representation}\\
\hline
$2^t-4 = 2^{11}-4$&$1$&$1$&$1$&$1$&$1$&$1$&$1$&$1$&$1$&$1$&$0$&$0$\\ 
$u = 2^5$&$0$&$0$&$0$&$0$&$0$&$0$&$1$&$0$&$0$&$0$&$0$&$0$\\ 
$v = 2^9+2^8-2^6$&$0$&$0$&$1$&$0$&$1$&$1$&$0$&$0$&$0$&$0$&$0$&$0$\\ 
$w = 2^{10}+2^8+2^5-4$&$0$&$1$&$0$&$1$&$0$&$0$&$0$&$1$&$1$&$1$&$0$&$0$\\ 
\hline
\end{tabular}
\end{table}
}
From Table \ref{tab3}, it is clear that $\binom{2^t-4}{u, v, w} = 1 \mod (2)$ (by Corollary \ref{lucascor}) for these $t-s$ values of $u$, i.e., for $u=0$ and $u=2^i$, with $2\leq i \leq t-s$. 
\item [Subcase 2.2.] ${\wt(u)> 1:}$ Let $u= \sum_{j=i_1-k+1}^{i_1} 2^j+ \sum_{j=0}^{i_2}u_j2^j$, where $k \geq 1, i_1 \leq t-s+1, i_2 \leq i_1-k-1, u_j \in \{0, 1\} ~\text{for}~0 \leq j \leq i_2$, and if $k=1$ then at least one $u_j$ is non-zero. Next, we consider the following subsubcases.
\begin{description}
\item [Subsubcase 2.2.1.] ${i_1 \leq t-s:}$ Note that $v= 2^{t-s+2}+2^{t-s+1}-2u$, i.e., $v = 2^{t-s+1}+(2^{t-s+2}-2u)$. Now, from the analysis of the Subsubcase 1.1.1 ($\wt(u)> 1$) of Case 1, we have that both $u$ and $2^{t-s+2}-2u$ has $1$ in the $r$-th bit position for some $r \leq i_1 \leq t-s$. Hence, both $u$ and $2^{t-s+2}+2^{t-s+1}-2u$ has $1$ in the same $r$-th bit position. So, following Corollary \ref{lucascor}, $\binom{2^t-4}{u, v, w} = 0 \mod (2)$.
\item [Subsubcase 2.2.2.]${i_1 = t-s+1:}$ We consider the following two possibilities.
\begin{enumerate}[label=(\textbf{\itshape\roman*})]
\item{$\wt(u - 2^{t-s+1}) \geq 2$:} Therefore, $v = 2^{t-s+2}+2^{t-s+1}-2u = 2^{t-s+1} - 2u'$, where $u' = u - 2^{t-s+1}$, and $\wt(u') \geq 2$.\linebreak Let $u' = 2^{j}+\sum_{m =0}^{j-1}u_{m}2^{m}$, $u_m \in \{0, 1\}$, and at least one $u_m = 1$. It is clear that $j \leq t-s-1$, for otherwise, $v < 0$. But, this is equivalent to the Subsubcase 1.1.1 ($\wt(u) > 1$) of Case 1, which implies $u'$ and $v$ has $1$ in the same $r$-th bit position for some $r \leq t-s-1$. This, in turn, implies $u$ and $v$ has $1$ in the same (as the previous) $r$-th bit position. Hence, we get that $\binom{2^t-4}{u,v,w} = 0 \mod (2)$ (by Corollary \ref{lucascor}) in this case as well.
\item{$\wt(u - 2^{t-s+1}) = 1$:} Let us assume that $u = 2^{t-s+1}+2^j$, where $2 \leq j \leq t-s$. In Table \ref{tab5}, we consider bit patterns of $u, v, w$ for this case.
\begin{table}[htb]
\caption{Bit patterns for Subsubcase 2.2.2-(ii)}
\label{tab5}
\hspace{1.5cm}
\begin{tabular}{|c|c|c|}
\hline
Values & \multicolumn{2}{c|}{Bit positions with $1$}\\
\hline
\multirow{3}{*}{
\pbox{20 cm}{
$u =2^{t-s+1}+2^j$,\\ $v = 2^{t-s+1} - 2^{j+1}$,\\ $w = 2^t - 2^{t-s+2}+2^j -4$\\
  $(2 \leq j \leq t-s)$ }} & $u$ & $\{j, t-s+1\}$\\ 
\cline{2-3}
& $v$ & $\{r \vert j+1 \leq r \leq t-s\}$\\
\cline{2-3}
& $w$ & $\begin{aligned}\big\{r \vert & 2 \leq r\leq j-1,\\
& t-s+2 \leq r \leq t-1\big\} \end{aligned}$\\
\hline
\end{tabular}
\end{table}
Bit patterns from Table \ref{tab5}, together with Corollary \ref{lucascor}, imply that for these $t-s-1$ values of $j$, $\binom{2^t-4}{u,v,w} = 1 \mod (2)$.
\end{enumerate} 
\end{description}
\end{description} 
So, considering the above possibilities, we conclude that the coefficient of the term with exponent $4(2^t-1)$ in the expansion of $f^{2^t-4}$ is 
\begin{eqnarray}
\label{2k4condn}
\alpha^{2^t-2^{t-s+2}-2^{t-s+1}-4}(1+\sum\limits_{i=2}^{t-s} \alpha^{2^i})+\alpha^{2^t-2^{t-s+2}-4} \sum\limits_{i=2}^{t-s}\alpha^{2^i}
\end{eqnarray}
Now, by employing (\ref{2k3condn}) and simplifying we get that the above expression equals $\alpha^{2^t-2^{t-s+2}+2^{t-s+1}-4}$, which is non-zero. From earlier discussion, the coefficient of the term with exponent $(2^{s-1}+4)(2^t-1)$ in the expansion of $f^{2^t-4}$ is $0$. So, the coefficient of the term with exponent $2^t-1$ in the expansion of $f^{2^t-4} \mod x^{2^t}+x$ is clearly non-zero.\qed \par\pagebreak
\bigskip
\subsection{Proof of Theorem \ref{binneg}}
{\em Proof.} We apply Theorem \ref{permth1}, where we set $g(x) = x(x^{2^n -1 \over 2^t-1}+a)$. Since $r=1$ in this case, so condition (a) of Theorem \ref{permth1} is satisfied.  Next, we observe that condition (b) of Theorem \ref{permth1} is satisfied if and only if $a\in \F_{2^{2t}}^{*} \setminus \F_{2^{t}}^{*}$. So, $g(x)$ is a PP if and only if $a\in \F_{2^{2t}}^{*} \setminus \F_{2^{t}}^{*}$ and condition (c) of Theorem \ref{permth1} is satisfied.\par
Let $\ga$ be a primitive element of $\F_{2^n}$, then $\beta = \ga^{2^n-1 \over 2^t-1}$ is a primitive element of $\F_{2^t}$. So, the condition $g(\ga^i)^{2^n -1 \over 2^t-1} \neq g(\ga^j)^{2^n -1 \over 2^t-1}$ for all $0 \leq i < j < 2^t-1$ is equivalent to the condition 
\begin{align}
\beta^i(\beta^i+a)^{{2^n-1} \over {2^t-1}} \neq \beta^j(\beta^j+a)^{{2^n-1} \over {2^t-1}},  ~ \mbox{for} ~i \neq j.
\label{eqcondn}
\end{align}
However, (\ref{eqcondn}) is equivalent to the condition that $x(x+a)^{{2^n-1} \over {2^t-1}}$ is a PP of $\F_{2^t}$.
 Hence, $g(x)$ is a PP of $\F_{2^n}$ if and only if $a \in \F_{2^{2t}}^{*} \setminus \F_{2^{t}}^{*}$ and $x(x+a)^{{2^n-1} \over {2^t-1}}$ is a PP of $\F_{2^t}$. Let $\F_{2^{2t}} = \F_{2^t}(\ze)$, where $\ze$ is a root of the irreducible polynomial $x^2+x+\theta \in\F_{2^t}[x]$. So, we have 
\begin{eqnarray}
\label{conjugate}
\ze+\ze^{2^t}=1 
\end{eqnarray}
Now, $a \in \F_{2^{2t}}^{*} \setminus \F_{2^{t}}^{*}$ can be written as $a = b  + c \ze$, where $b,c \in \F_{2^t}, c \neq 0$. So,
\begin{align}
x(x+a)^{{2^n-1} \over {2^t-1}} 	&= x(x+b+ c\ze)^{{2^{2^st}-1} \over {2^t-1}}\notag\\
     	&= x(x+b+ c\ze)^{2^{(2^s-1)t}+2^{(2^s-2)t}+\cdots+1}\notag\\
	&= x(x+b + c \ze)^{2^{(2^s-1)t}}(x+b + c \ze)^{2^{(2^s-2)t}}\ldots(x+b + c \ze)\notag\\
	&= x\prod\limits_{i=1}^{2^{s-1}}(x+b+c\ze^{2^t})(x+ b+ c \ze), \label{expand}
	\end{align}
where the expression in the r.h.s. of (\ref{expand}) follows by noting that for odd $\ell$, $(x+ b+ c \ze)^{2^{\ell t}} =x+b+c\ze^{2^t}$, and for even $\ell$, $(x+ b+ c \ze)^{2^{\ell t}} =x+b+c\ze$. Next, using (\ref{conjugate}), and after some regular calculation we have
\begin{align}
x\prod\limits_{i=1}^{2^{s-1}}(x+b+c\ze^{2^t})(x+ b+ c \ze)&= x\prod\limits_{i=1}^{2^{s-1}}(x^2+cx+b^2+bc+c^2\ze+c^2\ze^2).\label{usingconjugate}
\end{align}
Further, by applying the transformation $x \mapsto cy$ on the r.h.s. of (\ref{usingconjugate}), we have that $x(x+a)^{{2^n-1} \over {2^t-1}}$ is a PP of $\F_{2^{t}}$ if and only if the polynomial $c^{2^s}\left(y^{2^s+1}+\\y^{2^{s-1}+1}+\left({b^2+bc+c^2\ze+c^2\ze^2 \over c^2}\right)^{2^{s-1}}y\right)$ is a PP of $\F_{2^{t}}$.\par
 Since $b,c,~\text{and}~\ze^2+\ze = \theta \in \F_{2^t}^{*},$ we have ${b^2+bc+c^2\ze+c^2\ze^2 \over c^2} \in \F_{2^t}^{*}$. Now, by applying Theorem \ref{trith} on the polynomial $y^{2^s+1}+y^{2^{s-1}+1}+\left({b^2+bc+c^2\ze+c^2\ze^2 \over c^2}\right)^{2^{s-1}}y$ we have $x(x+a)^{{2^n-1} \over {2^t-1}}$ is a PP of $\F_{2^t}$ if and only if \begin{inparaenum}[(\itshape i\upshape)]\item $t$ is odd,\item $s\in \{1,2\}$, and \item ${{b^2+bc+c^2\ze+c^2\ze^2 \over c^2}} =1$, i.e., $b^2+bc+c^2(1+\ze+\ze^2)=0$. \end{inparaenum}  \par
Now, in order to express condition (\itshape iii \upshape) in terms of $a$, we observe (suppressing some regular calculation) that
\begin{align}
a^{2^{t+1}}+ a^2+ a^{2^t+1} &= (b+c\ze)^{2^{t+1}}+(b+c\ze)^2+(b+c\ze)^{2^t+1}\notag\\
				&= b^2+bc+c^2(1+\ze+\ze^2),\label{equivalent}
\end{align}
where we  employ (\ref{conjugate}) to derive the r.h.s. of (\ref{equivalent}). Therefore, condition (\itshape iii \upshape) is equivalent to the condition
\begin{align}
a^{2^{t+1}}+ a^2+ a^{2^t+1} = 0. \label{finalcond}
\end{align}
Finally, we make condition (\ref{finalcond}) succinct by noting that $a^{2^{t+1}}+ a^2+ a^{2^t+1} = a^2((a^{2^t-1})^2+a^{2^t-1}+1)$.
Therefore, $a^{2^{t+1}}+ a^2+ a^{2^t+1} = 0$ if and only if 
$a^{2^t-1}$ is a root of the equation ${x'}^2+x'+1 = 0$ 
in $\mathbb{F}_{2^{2t}}$, i.e., if and only if $a^{2^{t}-1} \in \{\omega, \omega^2\}$, which is true if and only if $a \in \om\F_{2^t}^{*}\bigcup \om^2\F_{2^t}^{*}$.
\qed  \par
\subsubsection*{Acknowledgement}
The authors thank anonymous reviewers for their helpful comments and corrections which improved the quality of this manuscript. The first author thanks Mr. Shashank Singh for commenting on an earlier draft.
\label{Bibliography}
\bibliographystyle{plain}  
\bibliography{dcc_revised}
\end{document}